\documentclass{amsart}
\usepackage{amsmath}
\usepackage{amscd}
\usepackage{amssymb}

\renewcommand{\div}{\operatorname{div}}

\newcommand{\Vol}{\operatorname{Vol}}

\newtheorem{theorem}{Theorem}[section]
\newtheorem{theorem/definition}{Theorem/Definition}[section]
\newtheorem{proposition}{Proposition}[section]

\newtheorem{corollary}{Corollary}[section]
\theoremstyle{remark}
\newtheorem{remark}{Remark}[section]

\theoremstyle{definition}
\newtheorem{definition}{Definition}[section]
\newtheorem{example}{Example}[section]

\begin{document}
\title
{Geometry of Complete Gradient Shrinking Ricci Solitons}
\author{Huai-Dong Cao}
\address{Department of Mathematics\\ Lehigh University\\
Bethlehem, PA 18015} \email{huc2@lehigh.edu}

%\begin{abstract}
%\end{abstract}

\dedicatory{Dedicated to Professor S.-T. Yau on the occasion of
his 60th birthday}

\maketitle
\date{}

%\footnotetext[1]{1991 {\em Mathematics Subject Classification}: }

\footnotetext[1]{Partially supported by NSF.}

The notion of {\it Ricci solitons} was introduced by Hamilton
\cite{Ha88} in mid 1980s. They are natural generalizations of
Einstein metrics. Ricci solitons also correspond to self-similar
solutions of Hamilton's Ricci flow \cite{Ha82}, and often arise as
limits of dilations of singularities in the Ricci flow. In this
paper, we will focus our attention on complete gradient shrinking
Ricci solitons and survey some of the recent progress, including
the classifications in dimension three, and asymptotic behavior of
potential functions as well as volume growths of geodesic balls in
higher dimensions.

\section{Gradient Shrinking Ricci Solitons}

Recall that a Riemannian metric $g_{ij}$ is {\it Einstein} if its
Ricci tensor $R_{ij}$ is a constant multiple of $g_{ij}$:
$R_{ij}=\rho g_{ij}$ for some constant $\rho$. A smooth
$n$-dimensional manifold $M^n$ with an Einstein metric $g_{ij}$ is
an {\it Einstein manifold}. Ricci solitons, introduced by Hamilton
\cite{Ha88}, are natural generalizations of Einstein metrics.

\begin{definition}
A complete Riemannian metric $g_{ij}$ on a smooth manifold $M^n$
is called a {\it gradient Ricci soliton} if there exists a smooth
function $f$ on $M^n$ such that the Ricci tensor $R_{ij}$ of the
metric $g_{ij}$ satisfies the equation
$$R_{ij}+\nabla_i\nabla_jf=\rho g_{ij}$$
for some constant $\rho$. For $\rho=0$ the Ricci soliton is {\it
steady}, for $\rho>0$ it is {\it shrinking} and for $\rho<0$ {\it
expanding}. The function $f$ is called a {\it potential function}
of the Ricci soliton.
\end{definition}

Note that when the potential function $f$ is a constant Ricci
soltions are simply Einstein metrics. In this paper, we will focus
our attention on complete gradient shrinking Ricci solitons, which
are possible Type I singularity models in the Ricci flow. We refer
the readers to \cite{Cao08b, Cao08a} and the references therein
for a quick overview and more information on singularity formation
of the Ricci flow on 3-manifolds and the role shrinking solitons
played.

Throughout the paper we will assume our gradient shrinking
solitons satisfy the equation
$$R_{ij}+\nabla_i\nabla_jf=\frac{1}{2} g_{ij}. \eqno(1.1)$$
This can be achieved by a suitable scaling of any shrinking
soliton metric $g_{ij}$.

To get a good feeling of how gradient shrinking solitons
correspond to self-similar ancient solutions of the Ricci flow, we
first observe that if $g_{ij}$ is an Einstein metric with positive
scalar curvature $R=n/2$, then $g_{ij}$ corresponds to a
homothetic shrinking solution ${g}_{ij}(t)$ to the Ricci flow
$$\frac{\partial g_{ij}(t)}{\partial t}=-2R_{ij}(t) $$
with $g_{ij}(0)=g_{ij}$. Indeed, if
$$R_{ij}=\frac{1}{2} g_{ij},$$  then
$${g}_{ij}(t)=(1-t)g_{ij} \eqno (1.2)$$ is a solution to the Ricci
flow which shrinks homothetically to a point as $t\rightarrow 1$.
Note that ${g}_{ij}(t)$ exists on the {\it ancient} time interval
$(-\infty, 1)$, hence an {\it ancient solution}, and the scalar
curvature ${R}(t)\to \infty$ like $1/(1-t)$ as $t\to 1$.
Similarly, a complete gradient shrinking Ricci soliton satisfying
equation (1.1) corresponds\footnote{It was observed by Z.-H. Zhang
\cite{Zh1} that $V=\nabla f$ is a complete vector field if the
soliton metric $g_{ij}$ is complete} to a self-similar solution
$\tilde{g}_{ij}(t)$ to the Ricci flow with
%\begin{align*}
$$\tilde{g}_{ij}(t):=(1-t)\varphi_t^*(g_{ij}),\qquad -\infty<t<1, \eqno(1.3)$$
%\end{align*}
where $\varphi_t$ is the 1-parameter family of diffeomorphisms
generated by $\nabla f/(1-t)$.

A theorem of Perelman states that given any non-flat
$\kappa$-noncollapsed ancient solution $g_{ij}(t)$ to the Ricci
flow with bounded and nonnegative curvature operator, the limit of
some suitable blow-back of the solution converges to a non-flat
gradient shrinking soliton (see Proposition 11.2 in \cite{P1}, or
Theorem 6.2.1 in \cite{CZ05}). Thus knowing the geometry of
gradient shrinking solitons helps us to understand the asymptotic
behavior of ancient solutions.

For compact shrinking Ricci solitons in low dimensions, we have

\begin{proposition} {\bf (Hamiton \cite{Ha88}\footnote{See alternative
proofs in (Proposition 5.21, \cite{CK}) or (Proposition 5.1.10,
\cite{CZ05}), and \cite{CLT}.} for $n=2$, Ivey
\cite{Iv}\footnote{See \cite{ELM} for alternative proofs} for
$n=3$)} In dimension $n\le 3$, there are no compact shrinking
Ricci solitons other than those of constant positive sectional
curvature.
\end{proposition}

\medskip
%\noindent{\bf Examples of gradient shrinking Ricci Solitons}

However, when $n\ge 4$, besides positive Einstein manifolds and
products of positive Einstein manifolds with Euclidean spaces,
there do exist non-Einstein compact gradient shrinking solitons.
Also, there exist complete noncompact non-flat shrinking solitons.
We list below the main examples.

\begin{example} {\bf Quotients of round n-sphere $\mathbb{S}^n/{\Gamma}$}

As we pointed out before, any positive Einstein manifold is a
(compact) shrinking Ricci soliton. In particular, the round
n-sphere $\mathbb{S}^n$, or any of its metric quotient
$\mathbb{S}^n/{\Gamma}$ by a finite group $\Gamma$, is a compact
shrinking soliton.
\end{example}

\begin{example} {\bf The Gaussian shrinker}

It is easy to check that the flat Euclidean space $(\Bbb R^n,
\delta_{ij})$ is a gradient shrinker with potential function
$f=|x|^2/4$:
$$\partial_i\partial_j f=\frac{1} {2}\delta_{ij}.$$
$(\Bbb R^n, \delta_{ij},|x|^2/4)$ is called the Gaussian shrinking
soliton.

\end{example}

\begin{example} {\bf The round cylinder $\mathbb{S}^{k}\times \mathbb{R}^{n-k}$}

The product of a positive Einstein manifold with the Euclidean
Gaussian shrinker is a (noncompact) shrinking soliton. In
particular, the round cylinder $\mathbb{S}^{k}\times
\mathbb{R}^{n-k}$ is a noncompact shrinking soliton.
\end{example}

\begin{example} {\bf Compact rotationally symmetric K\"ahler shrinkers}

For real dimension $4$, the first example of a compact shrinking
soliton was constructed in early 90's by Koiso \cite{Ko} and
independently by the author \cite{Cao94}\footnote{The construction
of the author was carried out in 1991 at Columbia University. When
he told his construction to S. Bando that year in New York, he
also learned the work of Koiso from Bando.} on compact complex
surface $\Bbb
 CP^2\#(-\Bbb CP^2)$, the below-up of the complex projective plane at one point.
Here $(-\Bbb CP^2)$ denotes the complex projective space with the
opposite orientation. This is a gradient shrinking K\"ahler-Ricci
soliton, which has $U(2)$ symmetry and positive Ricci curvature.
More generally, Koiso and the author found $U(n)$-invariant
K\"ahler-Ricci solitons on certain twisted projective line bundle
over $\Bbb CP^{n-1}$ for $n\geq 2$.

\end{example}

%\begin{remark} If a compact K\"ahler manifold $M$ admits a K\"ahler
%shrinker then $M$ is Fano (i.e., the first Chern class $c_1(M)$ of
%$M$ is positive), and the Futaki-invariant \cite{Fu} is nonzero.
%\end{remark}

\begin{example} {\bf Compact toric K\"ahler shrinkers}

In \cite{WZ}, Wang-Zhu found a gradient K\"ahler-Ricci soliton on
$\Bbb CP^2\#2(-\Bbb CP^2)$ which has $U(1)\times U(1)$ symmetry.
More generally, they proved the existence of gradient shrinking
K\"ahler solitons on all Fano toric varieties of complex dimension
$n\geq 2$ with non-vanishing Futaki invariant.
\end{example}

\begin{example} {\bf Noncompact gradient K\"ahler shrinkers}

Feldman-Ilmanen-Knopf \cite{FIK} found the first complete
noncompact $U(n)$-invariant gradient shrinking K\"ahler solitons
on certain twisted complex line bundles over $\Bbb CP^{n-1}$
($n\ge 2$) which are cone-like at infinity and have quadratic
curvature decay. Moreover, the shrinker has positive scalar
curvature but the Ricci curvature changes sign. The simplest such
example is defined on the tautological line bundle
$\mathcal{O}(-1)$ of $\Bbb CP^{1}$, the blow-up of $\Bbb C^{2}$ at
the origin.
\end{example}

\begin{remark}
Recently, Dancer-Wang \cite{DW1} produced new examples of
rotationally symmetric compact and complete noncompact gradient
K\"ahler shrinkers on bundles over the product of positive
K\"ahler-Einstein manifolds, generalizing those in Example 1.4 and
Example 1.6.
\end{remark}

Next we describe some useful results about gradient shrinking
solitons.

\begin{proposition} {\bf (Hamilton \cite{Ha95F})} Let $(M^n, g_{ij}, f)$
be a complete gradient shrinking soliton satisfying Eq. (1.1).
Then we have
$$\nabla_iR=2R_{ij}\nabla_jf, \eqno(1.4)$$ and
$$R+|\nabla f|^2-f=C_0 \eqno(1.5)$$ for some constant $C$. Here $R$
denotes the scalar curvature.
\end{proposition}

\begin{proof}
Let $(M^n, g_{ij}, f)$ be a complete gradient Ricci soliton
satisfying equation (1.1). Taking the covariant derivatives and
using the commutating formula for covariant derivatives, we obtain
$$\nabla_iR_{jk}-\nabla_jR_{ik}+R_{ijkl}\nabla_lf=0.$$
Taking the trace on $j$ and $k$, and using the contracted second
Bianchi,
%$$\nabla_j R_{ij}=\frac {1}{2}\nabla_i R,$$
we get
$$\nabla_iR=2R_{ij}\nabla_jf.$$
Thus
$$\nabla_i(R+|\nabla f|^2- f)=2(R_{ij}+\nabla_i\nabla_jf-\frac{1}{2} g_{ij})\nabla_jf=0.$$
Therefore $$R+|\nabla f|^2-f=C_0$$ for some constant $C_0$.
\end{proof}

Note that if we normalize $f$ by adding the constant $C_0$ to it,
then (1.5) becomes
$$R+|\nabla f|^2-f=0. \eqno(1.6)$$

%From now on, unless stated otherwise, we will make this
%normalization on the potential function $f$.

\begin{proposition} {\bf (Perelman \cite{P2}\footnote{This result was
essentially sketched by Perelman (see p.3 in \cite{P2}), and a
detailed argument was presented in \cite{CZ05} (see p.385-386 in
\cite{CZ05}).})} Let $(M^n, g_{ij}, f)$ be a complete noncompact
gradient shrinking Ricci soliton with bounded Ricci curvature
which satisfies (1.1) and the normalization (1.6). Let $r(x)=d(x,
x_0)$ denote the distance from $x$ to a fixed point $x_0\in M$.
Then there exist positive constants $C_1$, $C_2$, $c_1$ and $c_2$
such that, for $r(x)$ sufficiently large, $f$ satisfies the
following estimates:
$$\frac{1}{4}(r(x)-c_1)^2 \le f(x)\le C_1 (r(x)+c_2)^2; \eqno(1.7)$$
$$|\nabla f|(x)\le C_2 (r(x)+1).\eqno(1.8)$$
%Here $c_1>0$ is a positive constant depending only on $n$ and the
%geometry of $g_{ij}$ on the unit ball $B_{x_0}(1)$.
\end{proposition}

\begin{proof} By assumption, we know the Ricci curvature is bounded
by $|Rc|\le C$ for some positive constant $C>0$. From soliton
equation (1.1) and the Ricci lower bound $Rc\ge -Cg$, we
immediately have the upper estimate on Hessian of $f$:
$$\nabla_i\nabla_jf\le (C+1/2)g_{ij}, $$ from which the upper
estimate on $f$ in (1.7) follows. Now (1.8) follows from this
upper estimate $f(x)$, the scalar curvature lower bound $R\ge
-nC$, and (1.5) in Proposition 1.2.

To prove the lower estimate on $f$ in (1.7), consider any
minimizing normal geodesic $\gamma(s)$, $0\leq s \leq s_0$ for
some arbitrary large $s_0>0$ starting from $\gamma(0)=x_0$. Denote
by $X(s)=\dot\gamma(s)$ the unit tangent vector along $\gamma$.
Then, by the second variation of arc length, we have

$$\int_{0}^{s_0} \phi^{2}Rc (X, X) ds \leq (n-1) \int_{0}^{s_0}|\dot\phi(s)|^2 ds \eqno(1.9)$$
for every nonnegative function $\phi(s)$ defined on the interval
$[0, s_0]$. Now, following Hamilton \cite{Ha95F}, we choose $\phi
(s)$ by
$$ {\phi (s) =\left\{
       \begin{array}{lll}
  s, \ \ \qquad  s\in [0, 1],\\[4mm]
  1, \ \ \qquad  s\in [1, s_0-1], \\[4mm]
  s_0-s, \ \ s\in [s_0-1, s_0].
       \end{array}
    \right.}$$
Then
\begin{align*} \int_{0}^{s_0} Rc (X, X) ds
& = \int_{0}^{s_0} \phi^{2}Rc (X, X)ds + \int_{0}^{s_0}
(1-\phi^{2})Rc (X,
X)ds\\
& \leq (n-1) \int_{0}^{s_0}|\dot\phi(s)|^2 ds + \int_{0}^{s_0}
(1-\phi^{2})Rc (X,
X)ds\\
&\leq 2(n-1) + \max_{B_{x_0}(1)} |Rc| + \max_{B_{\gamma(s_0)}(1)}
|Rc|.
\end{align*}

On the other hand, by (1.1), we have
$$
\nabla_X\dot{f}= \nabla_X\nabla_Xf=\frac{1}{2}- Rc(X,X).
\eqno(1.10)
$$
Integrating (1.10) along $\gamma$ from $0$ to $s_0$, we get

$$\dot{f}(\gamma(s_0))-\dot{f}(\gamma(0))
=\frac{1}{2}s_{0}-\int_0^{s_0} Rc(X,X)ds. \eqno(1.11)$$

Since by assumption $g_{ij}$ has bounded Ricci curvature $|Rc|\le
C$, it follows that
\begin{align*}
\dot{f}(\gamma(s_0)) & \geq \frac{s_0}{2}+\dot{f}(\gamma(0)) -
2(n-1) - \max_{B_{x_0}(1)} |Rc| - \max_{B_{\gamma(s_0)}(1)} |Rc|\\
& \geq \frac{1}{2}s_{0}
-\dot{f}(\gamma(0))-2(n-1)-2C=\frac{1}{2}(s_{0}-c),
\end{align*}

and
$$ f(\gamma(s_0)) \geq \frac{1}{4}(s_{0}-c)^2-f(x_0)-\frac{c^2}{4}.$$
This completes the proof of Proposition 1.6.
\end{proof}

\begin{remark} When $Rc\ge 0$, (1.11) also implies that
$$f(x)\le \frac{1}{4} (r(x)+c_2)^2.\eqno (1.12)$$ On the other hand,
replacing $|Rc|\le C$ by $Rc\ge 0$, Ni \cite{Ni} showed that
$$ f(x)\ge \frac{1}{8} r^2(x)-C'_1. \eqno(1.13)$$
\end{remark}

\begin{remark} A precise asymptotic estimate on $f$,
without any curvature bound assumption, has been obtained more
recently by Cao-Zhou \cite{CZhou08} (see also Theorem 3.2 in this
paper).
\end{remark}

\section{Classification of 3-dimensional Gradient Shrinking Solitons}

In this section, we present the classification results on
3-dimensional complete shrinking Ricic solitons by Perelman
\cite{P2}, Ni-Wallach \cite{NW1}, and Cao-Chen-Zhu \cite{CCZ08}.

%First of all, we shall describe the classification result of
%Perelman \cite{P2} for 3-dimensional complete shrinking Ricic
%solitons which have {\it bounded and nonnegative sectional
%curvature} and are {\it $\kappa$-noncollapsed on all scales}. Then
%we will show, as done by Ni-Wallach \cite{NW1} and Cao-Chen-Zhu
%\cite{CCZ08}, how one can remove the assumptions on curvature
%bounds and $\kappa$-noncollaping.

In \cite{P2}, Perelman obtained the following classification
result of 3-dimensional complete shrinking solitons which is an
improvement of a result of Hamilton (Theorem 26.5, \cite{Ha95F}).

\begin{theorem} {\bf (Perelman \cite{P1})} Let $(M^3, g_{ij}, f)$ be
a complete nonflat 3-dimensional gradient shrinking Ricci soliton.
Suppose $g_{ij}$ has bounded and nonnegative sectional curvature
$0\le Rm \le C$ and is $\kappa$-noncollapsed on all scales for
some $\kappa>0$. Then $(M^3, g_{ij})$ is one of the following:

(i) the round three-sphere $\mathbb{S}^3$, or one of its metric
quotients;

(ii) the round infinite cylinder $\mathbb{S}^2\times \mathbb{R}$,
or one of its $\mathbb{Z}_2$ quotients.
\end{theorem}

We remark that, by the strong maximum principle of Hamilton
\cite{Ha86}, a complete shrinking Ricci soliton with nonnegative
curvature operator either has strictly positive curvature operator
everywhere or its universal cover splits as a product $N\times
\mathbb R^{k}$, where $k\ge 1$ and $N$ is a shrinking soliton with
positive curvature operator. On the other hand, we know that
compact shrinking solitons with positive curvature operator are
isometric to finite quotients of round spheres, thanks to the
works of Hamilton \cite{Ha82, Ha86} (for $n=3, 4$) and
B\"ohm-Wilking \cite{BW} (for $n\ge 5$).

In view of the proceeding remark, to prove Theorem 2.1 it suffices
to rule out the existence of a complete noncompact nonflat
$\kappa$-noncollapsed 3-dimensional gradient shrinking Ricci
soliton with bounded and positive sectional curvature.

\begin{proposition} {\bf (Perelman \cite{P1})} There is no
three-dimensional complete noncompact $\kappa$-noncollapsed
gradient shrinking soliton with bounded and positive sectional
curvature.
\end{proposition}

\medskip
\noindent {\it Sketch of the Proof}.\  By contradiction, suppose
there is a 3-dimensional complete noncompact
$\kappa$-non\-collapsed gradient shrinking soliton $(M^3, g_{ij},
f)$ with bounded and positive sectional curvature. %$g_{ij}(t)$,
%$-\infty<t<1$, be the corresponding self-similar solution as given in (1.).
Consider any minimizing normal geodesic $\gamma(s)$ starting from
$\gamma(0)=x_0$. As we saw in the proof of Proposition 1.3, for
$s$ sufficiently large we have
$$
\left|\nabla f\cdot\dot{\gamma}(s)-\frac{s}{2}\right|\le C,
$$
and
$$
\left|f(\gamma(s))-\frac{s^2}{4}\right|\le C\,\cdot(s+1).
$$
In particular, $f$ has no critical points outside some large
geodesic ball $B_{x_0}(s_0)$.

Now by (1.4) in Proposition 1.2 and the assumption of $Rc>0$, we
have

$$\nabla R\cdot \nabla f=2Rc (\nabla f, \nabla f)>0 \eqno(2.1)$$
at all points $x$ with $d(x,x_0)\ge s_0$. So outside
$B_{x_0}(s_0)$, the scalar curvature $R$ is strictly increasing
along the gradient curves of $f$. Hence

$$\bar{R}=\limsup_{d(x, x_0)\rightarrow+\infty}R(x)>0.\eqno(2.2)$$

\noindent {\bf Claim}: $ \bar{R}= 1$.
\medskip

First of all, if we choose a sequence of points $\{x_k\}$ so that
$R(x_k)\rightarrow\bar{R}$, then it follows from the noncollapsing
assumption and the compactness theorem that a subsequence of
$(M^3, g_{ij}, x_k)$ converges to some limit $(\bar{M}^3,
\bar{g}_{ij}, \bar{x})$, which splits off a line and hence is a
round cylinder with scalar curvature $\bar{R}$, corresponding to a
shrinking soliton defined (at least) on the ancient time interval
$(-\infty, 1)$. Thus we conclude that $\bar{R}\le 1$. In
particular, the scalar curvature of our original gradient
shrinking soliton satisfies the upper bound $$ R(x)< 1
\eqno(2.3)$$ outside $B_{x_0}(s_0)$.

To see that in fact $\bar{R}=1$, we consider the level surfaces of
$f$ $$\Sigma_s=\{x\in M: f(x)=s\}.$$ Note that the second
fundamental form of level surface $\Sigma_s$ is given by
$$h_{ij} = \nabla_i\nabla_j f/|\nabla f|,\qquad i,j=1,2,$$ where
$\{e_1,e_2\}$ is any orthonormal basis tangent to $\Sigma_s$. Now,
for 3-manifolds, the positivity of sectional curvature is
equivalent to $Rg_{ij}\ge 2R_{ij}$. If we choose an orthonormal
basis $\{e_1,e_2\}$ such that $Rc (e_1, e_2)=0$, then by (1.1) we
have
$$
\nabla_{e_i}\nabla_{e_j}f =\frac{1}{2}\delta_{ij}-Rc(e_i,e_j)
\ge(\frac{1}{2}-\frac{R}{2})\delta_{ij},\qquad i=1,2. \eqno(2.4)
$$ In particular, $\Sigma_s$ is convex and it follows that
\begin{align*} \frac{d}{ds}{\rm Area}\,(\Sigma_s)&
=\int_{\Sigma_s}\frac {1}{|\nabla
f|}\div(\frac{\nabla f}{|\nabla f|})\\
& \geq \int_{\Sigma_s}\frac{1}{|\nabla f|^2}(1-R)\\
& > \frac{1-\bar{R}}{2s}{\rm Area}\,(\Sigma_s)\ge 0
\end{align*}
for $s\ge s_0$. This implies that $\rm {Area}\,(\Sigma_s)$ is
strictly increasing as $s$ increases, and
$$
\log {\rm Area}\,(\Sigma_s)>(1-\bar{R})\log \sqrt{s}-C,
$$
for some constant $C$,  and $s\ge s_0$. But Area\,($\Sigma_s$) is
uniformly bounded from above by the area of the round sphere with
scalar curvature one.  Thus we conclude that $\bar{R}=1$, and
$$
{\rm Area}\,(\Sigma_s)<8\pi  \eqno(2.5)
$$ for $s$ large enough. This proves the Claim.

Now, denote by $\nu$ the unit normal to the level surface. By
using the Gauss equation and the soliton equation (1.1), the
intrinsic curvature $K$ of the level surface $\Sigma_s$ can be
computed as
\begin{align} K & =R_{1212}+\det(h_{ij})\\
& =\frac{R}{2}-Rc(\nu,\nu)+\frac{\det(\nabla^2 f)}{|\nabla f|^2}\\
& \le \frac{R}{2}-Rc(\nu,\nu)+\frac{(1-R+Rc(\nu,\nu))^2}{4|\nabla
f|^2} <\frac{1}{2}
\end{align}
for $s$ sufficiently large. But, this together with (2.5) lead to
a contradiction to the Gauss-Bonnet formula.

\medskip

\begin{corollary} The only three-dimensional complete noncompact
$\kappa$-noncollapsed gradient shrinking soliton with bounded and
nonnegative sectional curvature are either $\mathbb{R}^3$ or
quotients of $\mathbb{S}^2\times \mathbb{R}$.
\end{corollary}

Note that Perelman's proof of Proposition 2.1 relies upon the
Gauss-Bonnet formula, which cannot be used in higher dimensions.
In the past a few years, there have been considerable efforts to
improve or generalize the above results of Perelman. Ni-Wallach
\cite{NW1} and Naber \cite{Na} were able to drop the assumption on
$\kappa$-noncollapsing condition and replace nonnegative sectional
curvature assumption by that of nonnegative Ricci curvature. In
addition, Ni-Wallach \cite{NW1} allows the curvature $|Rm|$ to
grow as fast as $e^{ar(x)}$, where $r(x)$ is the distance function
to some origin $x_0$ and $a>0$ is a suitable small positive
constant.

\begin{theorem} {\bf (Ni-Wallach
\cite{NW1})} Any 3-dimensional complete noncompact non-flat
gradient shrinking soliton with nonnegative Ricci curvature
$Rc\geq 0$ and curvature bound $|Rm|(x)\leq Ce^{ar(x)}$ is a
quotient of the round sphere $\mathbb{S}^3$ or round cylinder
$\mathbb{S}^2\times \mathbb{R}$.
\end{theorem}

\medskip

\noindent{\it Sketch of the Proof}. \ It suffices to show that if
$(M^3, g_{ij}, f)$ is a complete gradient shrinking soliton with
positive Ricci curvature $Rc> 0$ and curvature bound $|Rm|(x)\leq
Ce^{ar(x)}$, then $(M^3, g_{ij})$ is a finite quotient of
$\mathbb{S}^3$. The basic ingredient of Ni-Wallach's proof is to
use the identity
$$\Delta (\frac{|Rc|^2}{R^2})=\nabla
(\frac{|Rc|^2}{R^2})\cdot\nabla f +\frac{2}{R^4}|R\nabla Rc
-\nabla R Rc|^2 - \frac{2}{R}\nabla
(\frac{|Rc|^2}{R^2})\cdot\nabla R + \frac{4P}{R^3} \eqno(2.7)$$
satisfied by the soltion metric, where
$$P=\frac{1}{2} ((\lambda+\mu-\nu)^2(\lambda-\mu)^2 +(\mu+\nu-\lambda)^2(\mu-\nu)^2
+(\nu+\lambda-\mu)^2(\nu-\lambda)^2),\eqno (2.8)$$ and $\lambda\ge
\mu\ge \nu$ are the eigenvalues of $Rc$. We remark that Hamilton
\cite{Ha82} showed that for any solution $g_{ij}(t)$ to the Ricci
flow on 3-manifolds there holds the evolution equation
$$ \frac{\partial}{\partial t} (\frac{|Rc|^2}{R^2})=\Delta
(\frac{|Rc|^2}{R^2}) -\frac{2}{R^4}|R\nabla Rc -\nabla R Rc|^2 +
\frac{2}{R}\nabla (\frac{|Rc|^2}{R^2})\cdot\nabla R -
\frac{4P}{R^3},$$ of which (2.7) is simply a special case when
$g_{ij}(t)$ is given by the self-similar solution (1.3)
corresponding to the shrinking soliton metric $g_{ij}$.

Now, by multiplying $|Rc|^2e^{-f}$ to (2.7) and integration by
parts, Ni-Wallach \cite{NW1} deduced that
$$0=\int_M (|\nabla
(\frac{|Rc|^2}{R^2})|^2R^2+\frac{2|Rc|^2}{R^4}|R\nabla Rc -\nabla
R Rc|^2+\frac{4P}{R^3}|Rc|^2)e^{-f}.$$ Thus: (i)
$\frac{|Rc|^2}{R^2}=constant$; (ii) $R\nabla Rc -\nabla R Rc=0$;
and (iii) $P=0$, provided the integration by parts is legitimate.
Moreover, it is clear from (2.8) that $Rc>0$ and $P=0$ imply
$\lambda=\mu=\nu$. Thus $R_{ij}=\frac{R}{3}g_{ij}$, which in turn
implies that $R$ is a (positive) constant and $(M^3, g_{ij})$ is a
space form.

Finally, based on the quadratic growth of $f$ in (1.13),
Ni-Wallach argued that the integration by parts can be justified
when the curvature bound $|Rm|(x)\leq Ce^{ar(x)}$ is satisfied.

Subsequently, B.-L. Chen, X.-P. Zhu and the author \cite{CCZ08}
observed that one can actually remove all the curvature bound
assumptions in Theorem 2.2.

\begin{theorem} {\bf (Cao-Chen-Zhu \cite{CCZ08})}
Let $(M^3, g_{ij}, f)$ be a $3$-dimensional complete non-flat
shrinking gradient soliton. Then $(M^3, g_{ij})$ is a quotient of
the round sphere $\mathbb{S}^3$ or round cylinder
$\mathbb{S}^2\times \mathbb{R}.$
\end{theorem}

\medskip

\noindent{\it Sketch of the Proof}. It was shown by B.-L. Chen
\cite{BChen} that every complete 3-dimensional ancient solution to
the Ricci flow, without assuming bounded curvature, has nonnegative sectional curvature $$Rm\ge 0.
\eqno(2.9)$$ In particular this means every complete noncompact
3-dimensional gradient shrinking soliton has $Rm\ge 0$.  Since we
now have $R\ge 0$, it follows from (1.6) that
$$0\leq |\nabla f|^2\leq f$$ Hence, $$|\nabla \sqrt{f}|\leq \frac{1}{2} \eqno(2.10)$$
whenever $f>0$. Thus $\sqrt{f}$ is an Lipschitz function, and it
follows that\footnote{The estimate here is more precise than that
given in \cite{CCZ08}, see also Lemma 2.3 in \cite{CZhou08}.}
$$\sqrt{f(x)}\leq \frac{1}{2} r(x) + \sqrt{f(x_0)},$$ where $x_0\in M^n$
is a fixed base point and $r(x)=d(x, x_0)$. Thus,
$$f(x)\leq \frac{1}{4} (r(x) +2\sqrt{f(x_0)})^{2}, \eqno(2.11)$$
$$|\nabla f|(x)\le \frac{1}{2}r(x)+\sqrt{f(x_0)},\eqno(2.12)$$ and
$$R(x)\le \frac{1}{4}(r(x)+2\sqrt{f(x_0)})^2. \eqno(2.13)$$
Consequently, $$|Rm|(x)\le  C (r^{2}(x) +1). \eqno(2.14)$$
Therefore, Theorem 2.3 follows from (2.9), (2.14) and Theorem 2.2.

\medskip
\begin{remark}
For $n=4$, Ni and Wallach \cite{NW2} showed that {\it any
$4$-dimensional complete gradient shrinking soliton with
nonnegative curvature operator and positive isotropic curvature,
satisfying certain additional assumptions, is a quotient of either
$\mathbb{S}^4$ or $\mathbb{S}^3\times \mathbb{R}$}. Based in part
on this result of Ni-Wallach, Naber \cite {Na} proved

\begin{theorem} {\bf (Naber \cite{Na})}
A 4-dimensional non-flat complete noncompact shrinking Ricci
soliton with bounded and nonnegative curvature operator is
isometric to a finite quotient of $\mathbb{S}^3\times \mathbb{R}$
or $\mathbb{S}^2\times \mathbb{R}^2$.
\end{theorem}

For $n\ge 4$, various results on classification of gradient
shrinking solitons with vanishing Weyl curvature tensor have been
obtained. In the compact case, see Eminenti-La Nave-Mantegazza
\cite{ELM}, while in the noncompact case, see Ni-Wallach
\cite{NW1}, Petersen-Wylie \cite{PW3}, X. Cao-Wang \cite{CW08},
and Z.-H. Zhang \cite{Zh2}. In particular, we have

\begin{proposition} {\bf (Ni-Wallach \cite{NW1})} Let $(M^{n},
g)$ be a complete, locally conformally flat gradient shrinking
soliton with nonnegative Ricci curvature. Assume that
$$|Rm|(x) \le e^{a(r(x) + 1)}$$ for some constant $a> 0$, where
$r(x)$ is the distance function to some origin. Then its universal
cover is  $\mathbb{R}^n$, $\mathbb{S}^n$ or
$\mathbb{S}^{n-1}\times \mathbb{R}$.
\end{proposition}

\begin{proposition} {\bf (Petersen-Wylie \cite{PW3})}
Let $(M^{n}, g)$ be a complete gradient shrinking Ricci soliton
with potential function $f$. Assume the Weyl tensor $W = 0$ and
$$\int_{M} |Rc|^{2} e^{-f} <\infty,$$ then $(M^{n}, g)$ is a finite
quotient of $\mathbb{R}^n$, $\mathbb{S}^n$ or
$\mathbb{S}^{n-1}\times \mathbb{R}$.
\end{proposition}

In \cite{Zh2}, Z.-H. Zhang showed that an $n$-dimensional ($n\ge
4$) complete gradient shrinking soliton with vanishing Weyl tensor
necessarily has nonnegative curvature operator. Thus, combining with Proposition 2.2, we have

\begin{proposition} {\bf (Z.-H. Zhang \cite{Zh2})}
Any $n$-dimensional ($n\ge 4$) complete non-flat gradient
shrinking soliton with vanishing Weyl tensor must be a finite
quotient of $\mathbb{S}^n$ or $\mathbb{S}^{n-1}\times \mathbb{R}$.
\end{proposition}
\end{remark}

\begin{remark}
Most recently, we have learned that Munteanu-Sesum \cite{MS09} have shown that
$\int_{M} |Rc|^{2} e^{-f} <\infty$ for any complete gradient shrinking soliton. Thus, in view of 
Proposition 2.3, this gives an alternative proof of Proposition 2.4. We refer the readers to 
\cite{MS09} for more information.  
\end{remark}

\section{Geometry of complete gradient solitons}

In the previous section, we described the classification results
on complete gradient shrinking solitons in dimensions $n=3$ and
$n=4$, as well as in the locally conformal flat case when $n\ge
4$. In all three cases, the proofs depend on the crucial fact
that the gradient shrinking soliton under the consideration has
nonnegative curvature operator. We also saw the condition of
nonnegative curvature operator is automatically satisfied when
$n=3$ or when the shrinker is locally conformally flat. However,
when $n\ge 4$ one cannot expect that in general a complete
gradient shrinking soliton has nonnegative sectional curvature, or
even nonnegative Ricci curvature. For example, the noncompact
gradient K\"ahler shrinkers of Feldman-Ilmanen-Knopf \cite{FIK}
does not have nonnegative Ricci curvature. While it is certainly
of considerable interest to classify complete gradient shrinking
solitons with nonnegative curvature operator or nonnegative Ricci
curvature in dimension $n\ge 5$, it is also important to
understand as much as possible the geometry of a general complete
noncompact shrinker $(M^n, g_{ij}, f)$ for $n\ge 4$. In this
section, we report some recent progress in this direction. The
results we are going to describe are mainly concerned with the
asymptotic behavior of potential functions, and volume growth
rates of geodesic balls.

First of all, we shall need the logarithmic Sobolev inequality
proved recently by Carrillo-Ni \cite{CN} for complete noncompact
gradient shrinking solitons. To state their result, instead of
normalizing the potential function $f$ by (1.6), we will use a
different normalization by requiring

$$\int_M e^{-f} =(4\pi)^{\frac{n}{2}}.\eqno(3.1)$$

\begin{remark}
According to \cite{M, WW} or \cite{CZhou08},  $\int_M e^{-f}$ is
finite on a complete gradient shrinking soliton $(M^n, g_{ij},
f)$. Hence normalization (3.1) is always possible.
\end{remark}

Now under the normalization (3.1) on $f$, we have
$$ R+|\nabla f|^2 -f=\mu_0\eqno(3.2)$$ for some (not necessarily positive)
constant $\mu_0$.

\begin{proposition} {\bf (Carrillo-Ni \cite{CN})} Let $(M^n, g_{ij}, f)$ be
a complete noncompact gradient shrinking soliton satisfying (1.1).
Assume that the scalar curvature is bounded below by a negative
constant. Then, for any function $u\in C^{\infty}_0(M)$ with
$\int_M u^2=(4\pi)^{n/2}$, one has
$$(4\pi)^{-\frac{n}{2}}\int_M (Ru^2+4|\nabla u|^2-u^2\log
u^2-nu^2)\ge \mu_0, \eqno(3.3)$$ where $\mu_0$ is the same
constant as given in (3.2).
\end{proposition}

We shall also need the following very useful fact\footnote{This is
a special case of Proposition 5.5 stated in \cite{Cao08b} which
was essentially proved in \cite{BChen} .} due to B.-L. Chen \cite{BChen}.

\begin{proposition}{\bf (B.-L. Chen \cite{BChen})} Let $(M^n, g_{ij}, f)$
be a complete shrinking Ricci soliton. Then the scalar curvature
$R$ of $g_{ij}$ is nonnegative: $$R\ge 0.$$
\end{proposition}

As a consequence of Propositions 3.1 and 3.2, one can deduce the
following Perelman-type non-collapsing result, which is a slight
improvement of Corollary 4.1 in \cite{CN}, for gradient shrinking
solitons by using a similar argument as in the proof of Theorem
3.3.3 in \cite{CZ05}.

\begin{corollary}  Let $(M^n, g_{ij}, f)$
be a complete noncompact gradient shrinking soliton satisfying
(1.1).  Then there exists a positive constant $\kappa=\kappa
(\mu_0)>0$ such that whenever $r\le 1$ and $R\le \frac{C}{r^2}$ on
a geodesic ball $B_r\subset M$, one has $\Vol(B_r) \ge \kappa
r^n.$
\end{corollary}

Now we are ready to prove the following result obtained by X.-P.
Zhu and the author \cite{CZ08}.

\begin{theorem} {\bf (Cao-Zhu \cite{CZ08})}
Let $(M^n, g_{ij}, f)$ be a complete noncompact gradient shrinking
Ricci soliton. Then $(M^n, g_{ij})$ has infinite volume:
$$\Vol(M^n, g_{ij})=\infty.$$
 More precisely, there exists some positive constant $C_3>0$ such
that
$$ \Vol(B_{x_0}(r))\ge C_3 \ln\ln r$$
for $r>0$ sufficiently large.
\end{theorem}

\begin{proof} We are going to show that if $\Vol(M^n, g_{ij})< \infty$,
then we shall get a contradiction to the logarithmic Sobolev
inequality (3.3). The argument is similar in spirit to that of the
proof of Theorem 3.3.3 in \cite{CZ05}, and is an adoption of the
Perelman's proof on the uniform diameter estimate for the
normalized K\"ahler Ricci flow on Fano manifolds.
%(cf \cite{Sesum}).

Let $(M^n, g_{ij}, f)$ be any complete noncompact gradient
shrinking soliton satisfying equation (1.1) and the normalization
condition (3.1). Then, by Proposition 3.2, we have $R\ge 0$. Thus
the log Sobolev inequality (3.3) of Carrillo-Ni is valid for
$(M^n, g_{ij}, f)$.

Pick a base point $x_0\in M$ and denote by $r(x)=d(x, x_0)$ the
distance from $x$ to $x_0$. Also denote by $A(k_1, k_2)$ the
annulus region defined by
$$A(k_1, k_2)=\{x\in M: 2^{k_1}\le d(x, x_0)\le 2^{k_2}\},$$
and $$V(k_1, k_2)=\Vol (A(k_1, k_2)),$$ the volume of $A(k_1,
k_2)$.

Note that for each annulus $A(k, k+1)$, the scalar curvature is
bounded above by $R\le C2^{2k}$ for some uniform constant $C>0$.
Since $A(k, k+1)$ contains at least $2^{2k-1}$ balls $B_r$ of
radius $r=2^{-k}$ and each of these $B_r$ has $\Vol(B_r)\ge \kappa
(2^{-k})^n$ by Corollary 3.1,
$$ V(k, k+1)\ge \kappa 2^{2k-1}2^{-kn}. \eqno(3.4)$$

Now, suppose $$\Vol(M^n, g_{ij})< \infty.$$ Then, for each
$\epsilon
>0$, there exists a large constant $k_0>0$ such that if
$k_2>k_1>k_0$, then $$V(k_1, k_2)\le \epsilon. \eqno(3.5)$$  We
claim that we can choose $k_1$ and $k_2$ in such a way that
$V(k_1, k_2)$ also satisfies the following volume doubling
property:
$$V(k_1, k_2)\le 2^{4n}V(k_1+2, k_2-2).\eqno(3.6)$$

Indeed, we can choose a very large number $K>0$ and pick
$k_1\approx K/2, k_2\approx 3K/2$. Suppose (3.6) does not hold,
i.e.,
$$V(k_1, k_2)> 2^{4n}V(k_1+2, k_2-2).$$
Then we consider whether or not
$$V(k_1+2, k_2-2)\le 2^{4n}V(k_1+4, k_2-4).$$
If yes, then we are done. Otherwise we repeat the process. After
$j$ steps, we have
$$V(k_1, k_2)> 2^{4nj}V(k_1+2j, k_2-2j).$$

However, when $j\approx K/4$ and using (3.4), this implies that
$$\Vol(M^n)>V(k_1, k_2)\ge 2^{nK}V(K, K+1)\ge \kappa 2^{2K-1}.$$
But $\Vol(M^n)$ is suppose to be finite, so after finitely many
steps (3.6) must hold for a pair of large numbers $k_2>k_1$. Thus
we can choose $k_1=k_1(\epsilon)$ and $k_2\approx 3k_1$ so that
both (3.5) and (3.6) are valid.

Next by using the Co-Area formula and the fundamental theorem of
calculus, we can choose $r_1\in [2^{k_1}, 2^{k_1+1}]$ and $r_2\in
[2^{k_2}, 2^{k_2+1}]$ such that
$$\Vol(S_{r_1})\le \frac{V(k_1, k_2)}{2^{k_1}} \quad \mbox{and}
\quad \Vol(S_{r_2})\le \frac{V(k_1, k_2)}{2^{k_2}}, $$ where
$S_{r}$ denotes the geodesic sphere of radius $r$ centered at
$x_0$.

Then, by integration by parts and using (2.12),
\begin{align*} |\int_{A(r_1, r_2)}\Delta f|& \le
\int_{S_{r_1}}
|\nabla f| +\int_{S_{r_2}} |\nabla f|\\
& \le \frac{V(k_1, k_2)}{2^{k_1}}C2^{k_1+1} + \frac{V(k_1,
k_2)}{2^{k_2}}C2^{k_2+1}\\
& \le C V(k_1, k_2).
\end{align*}
Therefore, since $R+\Delta f=n/2$, it follows that
$$\int_{A(r_1, r_2)} R\le CV(k_1, k_2),\eqno (3.7)$$
for some universal positive constant $C>0$.

Now we can derive a contradiction to the log Sobolev inequality
(3.3). Pick a smooth cut-off function $0<\zeta(t)\leq 1$ defined
on the real line such that

$$ {\zeta (t) =\left\{
       \begin{array}{lll}
  1, \ \ \qquad  2^{k_1+2}\le t\le 2^{k_2-2},\\[4mm]
  0, \ \ \qquad \mbox{outside}\ [r_1, r_2],
       \end{array}
    \right.}$$
and $|\zeta'|\leq 1$ everywhere. Define
$$u=e^{L}\zeta({d(x,x_0)}),
$$
where the constant $L$ is chosen so that $$(4\pi)^{n/2}=\int_M
u^2=e^{2L}\int_{A(r_1, r_2)}\zeta^2.\eqno(3.8)$$ Then, by the log
Sobolev inequality (3.3), we have

\begin{align*} \mu_0 & \le (4\pi)^{-\frac{n}{2}}\int_M (Ru^2+4|\nabla u|^2-u^2\log
u^2-nu^2)\\
& = (4\pi)^{-\frac{n}{2}}e^{2L}\int_{A(r_1, r_2)}R \zeta^2
+4|\zeta'|^2-2\zeta^2\log\zeta-2\zeta^2\log L-2n\zeta^2\\
& = -2\log L-2n + (4\pi)^{-\frac{n}{2}}e^{2L}\int_{A(r_1, r_2)} R
\zeta^2 + 4|\zeta'|^2-2\zeta^2\log\zeta.\\
\end{align*}

Now we can estimate the two remaining integral as follows: by
(3.7) and (3.6) we have
\begin{align*}
e^{2L}\int_{A(r_1, r_2)} R \zeta^2 & \le
Ce^{2L}V(k_1, k_2)\le Ce^{2L}2^{8n}V(k_1+2, k_2-2)\\
& \le C2^{8n}\int_{A(r_1, r_2)} u^2 \le
C2^{8n}(4\pi)^{\frac{n}{2}};
\end{align*}
on the other hand, using the elementary fact that $-t\log t\le e$
for $0\le t\le 1$, we obtain
\begin{align*}
e^{2L}\int_{A(r_1, r_2)}4|\zeta'|^2-2\zeta^2\log\zeta & \le
Ce^{2L}V(k_1, k_2)\\
& \le Ce^{2L}2^{8n}V(k_1+2, k_2-2)\\
& \le C2^{8n}\int_{A(r_1, r_2)} u^2 \le
C2^{8n}(4\pi)^{\frac{n}{2}}.
\end{align*}
Therefore,
$$ \mu_0\le -2\log L-2n-C,$$ for some universal positive constant $C>0$.
But, this is a contradiction, since by (3.5) and (3.8) we can make
$L$ arbitrary large by choosing $\epsilon>0$ arbitrary small.

Moreover, if one examines the above proof carefully, one can see
that in fact the geodesic balls $B_r(x_0)$ necessarily have at
least $\ln\ln r$ growth. This completes the proof of Theorem 3.1.
\end{proof}

\medskip
More recently, D. Zhou and the author \cite{CZhou08} obtained a
rather precise estimate on asymptotic behavior of the potential
function of a general complete gradient shrinking soliton.

\begin{theorem} {\bf (Cao-Zhou \cite{CZhou08})} Let $(M^n, g_{ij}, f)$
be a complete noncompact gradient shrinking Ricci soliton
satisfying (1.1) and the normalization (1.6). Then, the potential
function $f$ satisfies the estimates
$$\frac{1}{4} (r(x)-c_3)^2\leq f(x)\leq \frac{1}{4} (r(x)+c_4)^2.\eqno(3.9)$$
Here $r(x)=d(x_0, x)$ is the distance function from some fixed
point $x_0\in M$, $c_3$ and $c_4$ are positive constants depending
only on $n$ and the geometry of $g_{ij}$ on the unit ball
$B_{x_0}(1)$.
\end{theorem}

\begin{remark} In view of the Gaussian shrinker in Example 1.2, whose
potential function is $|x|^2/4$, the leading term
$\frac{1}{4}r^{2}(x)$ in (3.9) is optimal.
\end{remark}

\begin{remark} By Theorem 3.2, we can use $f$ to define a distance-like
function on $M^n$ by
$$\rho(x)=2\sqrt{f(x)} \eqno(3.10)$$ so that
$$r(x)-c\leq \rho(x)\leq r(x)+c \eqno(3.11)$$ with
$c=\max \{c_3, c_4\}>0$. Moreover, by (1.6) and Proposition 3.2,
we have $|\nabla f|^2\leq f$. Hence, whenever $f>0$,
$$|\nabla \rho|=\frac{|\nabla f|}{\sqrt{f}}\le 1.\eqno(3.12)$$
Note also that (1.6) and the upper bound on $f$ in (3.9) imply
that
$$R(x)\leq \frac{1}{4}(r(x)+2\sqrt{f(x_0)})^2. \eqno(3.13)$$
\end{remark}

Using Theorem 3.2 and the soliton equation (1.1), and working with
$\rho(x)$ as a distance function, an upper estimate on the volume
growth was derived in \cite{CZhou08}.

\begin{theorem} {\bf (Cao-Zhou \cite{CZhou08})} Let $(M^n, g_{ij}, f)$
be a complete noncompact gradient shrinking Ricci soliton. Then, there exists some positive constant $C_4>0$ such that
%$$(i) \quad \Vol(B_{p}(r))\geq C_1 r, \ \mbox{and}$$
$$ \Vol(B_{x_0}(r))\leq C_4 r^{n}$$
for $r>0$ sufficiently large.
\end{theorem}

\begin{remark}
  The noncompact K\"ahler shrinker of
Feldman-Ilmanen-Knopf \cite{FIK} in Example 1.6 has Euclidean
volume growth, with $Rc$ changing signs and $R$ decaying to zero.
This shows that the volume growth rate in Theorem 3.3 is optimal.
\end{remark}

\begin{remark}
Carrillo-Ni \cite {CN} proved that any non-flat gradient shrinking
soliton with {\sl nonnegative Ricci curvature $Rc\ge 0$} must have
zero asymptotic volume ratio, i.e., $\lim_{r\to\infty}
\Vol(B_{x_0}(r))/r^n=0.$
\end{remark}

In addition, the following result was also shown by Cao-Zhou
\cite{CZhou08}.

\begin{proposition} {\bf (Cao-Zhou \cite{CZhou08})} Let $(M^n, g_{ij}, f)$ be a
complete noncompact gradient shrinking Ricci soliton. Suppose the average
scalar curvature satisfies the upper bound
$$ \frac{1}{V(r)}\int_{D(r)} R\leq \delta$$
for some
constant $0<\delta < n/2$, and all sufficiently large r. Then, there exists some positive
constant $C_5>0$ such that
$$\Vol(B_{x_0}(r))\geq C_5 r^{n-2\delta}$$
for $r$ sufficiently large. Here $D(r)=\{x\in M^n: \rho(x)< r\}$, with $\rho(x)$ given by (3.10), and  $V(r)=\Vol(D(r))$.
\end{proposition}

\begin{remark}
On a complete noncompact Riemannian manifold $X^n$ with
nonnegative Ricci curvature, a theorem of Yau and Calabi (see
\cite{Yau}) asserts that the geodesic balls of $X$ have at least
linear growth, while the classical Bishop volume comparison
theorem says (cf. \cite{SY}) the geodesic balls of $X$ have at
most Euclidean growth. Theorem 3.1 and Theorem 3.3 were motivated
by these two well-known theorems respectively. It remains an
interesting problem to see if every complete noncompact gradient
shrinking soliton has at least linear volume growth. 
\end{remark}

\begin{remark}
On the topological side, Wylie \cite{Wy} showed that a complete
shrinking Ricci soliton has finite fundamental group. In the
compact case, this result was proved by\footnote{As pointed out in
\cite{WW}, the finiteness of fundamental group in the compact case
was shown in the earlier work of X. Li \cite{Li} by using
probabilistic methods.} Derdzinski \cite{De06}, and
Fern\'adez-L\'opez and Garc\'ia-R\'io \cite{FG} (see also a
different proof by Eminenti-La Nave-Mantegazza \cite{ELM}).
Moreover, Fang-Man-Zhang \cite{FMZ} proved that a complete
gradient shrinking Ricci soliton with bounded scalar curvature has
finite topological type.
\end{remark}

\medskip \noindent {\bf Acknowledgments}. I would like to dedicate
this paper to my teacher Professor S.-T. Yau, who taught me so
much and guided me for over a quarter of century. I would also
like to thank X.-P. Zhu for very helpful discussions.

\end{document}